\documentclass[12pt]{article}
\usepackage{amsmath,amsfonts}
\def\be{\begin{equation}} 
\def\ee{\end{equation}}
\newcommand{\beq}{\begin{eqnarray}}
\newcommand{\eeq}{\end{eqnarray}}
\newcommand{\nbeq}{\begin{eqnarray*}}
\newcommand{\neeq}{\end{eqnarray*}}

\begin{document}

\title{Characterizations of Logistic Distribution Through Order Statistics with Independent Exponential Shifts}
\author{M. Ahsanullah  \\
{\it Rider University, Lawrenceville, USA} \\ \\ George P. Yanev and Constantin Onica\\ {\it The University of Texas - Pan American, Edinburg,
USA}\\ \\  } \date{} \maketitle

\small{
\begin{tabular}{ll}

\noindent  M. Ahsanullah & George P. Yanev and Constantin Onica\\
 Department of Management Sciences & Department of Mathematics \\
  Rider University & University of Texas - Pan American \\
 2083 Lawrenceville Road & 1201 West University Drive  \\
Lawrenceville, NJ 08648, USA & Edinburg, TX 78539, USA \\
 E-mail: ahsan@rider.edu & E-mail: yanevgp@utpa.edu; onicac@utpa.edu
\end{tabular}

}

\newpage

\begin{abstract}
The paper presents some distributional properties of logistic order statistics subject to independent exponential one-sided and two-sided shifts.
Utilizing these properties, we extend several known results and obtain some
new characterizations of the logistic distribution.
 \\


\end{abstract}

\section{Discussion of Main Results}

The standard logistic cumulative distribution function (c.d.f.) is given by
\begin{equation}
\label{log_cdf} F_{L}(x)=\frac{e^x}{1+e^x}, \qquad -\infty <x<\infty.
\end{equation}
The density curve resembles that of normal distribution with heavier tails (higher kurtosis). The mean and variance are 0 and $\pi^2/3$, respectively. The logistic distribution has an important place in both probability theory and statistics as a model for population growth. It has been successfully applied to modeling in such diverse areas as demography, biology, epidemiology, environmental studies, psychology, marketing, etc.  Applications in economics include portfolio modeling Bowden \cite{Bo06}, approximation of the fill rate of inventory systems Zhang and Zhang \cite{ZZ07}, and Hubbert models of production trends of various resources Modis \cite{Mo98}, among others.

Characterizations of distributions is an active area of contemporary probability theory. They reveal intrinsic properties of distributions as well as connections, sometimes unexpected, between them. Within the abandon literature, there are surprisingly few results about the logistic distribution (see the discussion in  Galambos \cite{Ga92}  and Lin and Hu \cite{LH08}). George and Mudolkar \cite{GM81} were first to obtain characterizations involving order statistics.
Recently, Lin and Hu \cite{LH08} generalized their results and
derived new interesting ones. In this paper,
we extend  some of George and Mudolkar's \cite{GM81} and
Lin and Hu's \cite{LH08} findings, see Theorem 2 and Corollary below. Furthermore, we obtain new characterizations making use of distributional relationships among order statistics with independent exponential shifts.

Let $X_1, X_2,
\ldots, X_n $ be independent copies of a random variable $X$ with absolutely continuous (with respect to Lebesgue measure) c.d.f. $F(x)$. Denote the corresponding order statistics by $X_{1,n}\le X_{2,n}\le  \ldots\le X_{n,n}$. Let $\phi_{k}(t)$  and $F_k(x)$ (suppressing the dependence on $n$) be the characteristic function (ch.f.) and the c.d.f. of $X_{k,n}$ for $1\le k\le n$, respectively.
Further on $"\stackrel{d}{=}"$ stands for equality in distribution. Finally, let $E'_j$ and $E''_j$, $1\le j\le n$, be mutually independent standard exponential random variables, which are also  independent from $X_i$ for $1\le i\le n$.

In the first theorem below, we characterize the logistic distribution by distribution equalities between two order statistics with random exponential shifts. Three cases are considered depending on the relative positions of the order statistics: adjacent, two, and three spacings away.

\vspace{0.5cm}{\bf Theorem 1 }\ {\it Suppose $F(x)$ is an absolutely continuous (with respect to Lebesgue measure) c.d.f. with $F(0)=1/2$. Choose $r$ to be 1, 2, or 3 and let $k_i$ for $1\le i\le r$ be distinct integers in $[1,n-r]$.  Then $F$ is standard logistic if and only if

(i)  $t^r\phi_{k_i}(t)$ and $t^r\phi_{k_i+r}(t)$
are absolutely integrable for any $t$; and

(ii) the following $r$ equations hold for $k_i\in [1,n-r]$ where $1\le i\le r$
\be \label{lemmaA}
X_{k_i+r,n}-\sum_{j=k_i}^{k_i+r-1}\frac{E'_{j}}{j}\stackrel{d}{=}X_{k_i,n}+\sum_{j=k_i}^{k_i+r-1}\frac{E''_{j}}{n-j}.
\ee
}

{\bf Remark}\ It is worth mentioning here Corollary~2 of Weso{\l}owski and Ahsanullah \cite{WA04} where, in the context of order statistics with power contractions, the exponential distribution is characterized by $X_{k+1,n}\stackrel{d}{=}X_{k,n}+E'_1/(n-k)$ for $1\le k\le n-1$.

\vspace{0.5cm}Our second theorem
 extends the characterization
$X\stackrel{d}{=} X_{1,n} +\sum_{j=1}^{n-1}E'_j/j$ given in Lin and Hu \cite{LH08} by substituting $X_{1,n}$ with
$X_{k,n} -\sum_{j=1}^{k-1}E''_j/j$ for $1\le k\le n-1$.

\vspace{0.5cm}{\bf Theorem 2}\  {\it Let $F(x)$ be absolutely continuous (with respect to Lebesgue measure) and $F(0)=1/2$. Assume

(i) $G_X(x)=P(e^X\le x)$ is analytic and strictly increasing in $[0,\infty)$;

(ii) The derivative $G_X^{(k)}(x)$, $k\ge 1$, is strictly monotone in some interval $[0,t_k)$.

\noindent Choose an arbitrary but fixed $k\in \{1,2,\ldots n-1\}$. Then $F$ is standard logistic
iff
\be \label{goal1}
X\stackrel{d}{=} X_{k,n} +\sum_{j=1}^{n-k}\frac{E'_j}{j}-\sum_{j=1}^{k-1}\frac{E''_j}{j},
\ee
where as usual $\sum_{j=1}^0 (\cdot)=0$.}

\vspace{0.3cm} It is worth noticing the following corollary concerning the sample median $X_{k,2k-1}$. It generalizes George and Mudholkar's \cite{GM81}
$
X\stackrel{d}{=}X_{2,3}+E'_1-E'_2
$
characterization (see also Kotz et al. \cite{KKP01}, (2.7.37)).

\vspace{0.5cm}{\bf Corollary } {\it Let $k\ge 2$ and the assumptions in Theorem~2 hold. Then $F$ is standard logistic if and only if
\be \label{med_identity}
X\stackrel{d}{=} X_{k,2k-1} +\sum_{j=1}^{k-1}{La}_j,
\ee
where $La_j$  for $1\le j\le k-1$ are mutually independent Laplace random variables with density function $f_j(x)=je^{-j|x|}/2$ for $|x|<\infty$ and independent from $X_{k,2k-1}$.}

\vspace{0.3cm}In Section 2 we list some auxiliary results as well as derive two distributional relations among logistic order statistics. In the next two sections we consider the adjacent and non-adjacent cases of Theorem 1, respectively. Section 5 deals with the proof of Theorem 2 making use of the method of intensively monotone operators.


\section{Preliminaries}

In this section we gather some known formulas related to order statistics of absolutely continuous distributions, which will be used in  the next sections.  Let us start with the following characterization property of the logistic distribution (cf. Lin and Hu \cite{LH08}).
If $F$ is an absolutely continuous c.d.f. such that
\be \label{property}
F'(x)=F(x)(1-F(x),\qquad  x\in (l_F,r_F),
\ee
where $l_F=\inf\{x:F(x)>0\}$ and $r_F=\sup\{x:F(x)<1\}$,  then  $F$ is logistic given by $F(x)=\left(1+e^{-(x-\mu)}\right)^{-1}$ for some constant $\mu$ and $x\in R$.

We will also need some relations between the c.d.f.'s of $X$ and its order statistics. First, recall that (e.g., Ahsanullah and Nevzorov \cite{AN01}, Eqn. (1.1.2))
\be \label{os_density9}
F'_{k}(x)=k{n \choose k}F^{k-1}(x)(1-F(x))^{n-k}F'(x), \qquad 1\le k\le n.
\ee
Integrating (\ref{os_density9}) and iterating, we obtain the recurrence for $1\le r\le n-k$
\be \label{os_density3}
F_{k}(x) 
  =
 \sum_{j=k}^{k+r}{n \choose j}F^{j}(x)(1-F(x))^{n-j}+ F_{k+r+1}(x).
 \ee
Next, we recall some inversion formulas for characteristic functions. If $F_{j}(x)$ for $1\le j\le n$ is absolutely continuous, then
\[
F'_{j}(x)=\frac{1}{2\pi}\int_{-\infty}^\infty e^{-itx}\phi_{j}(t)dt, \quad F'_{j}(\infty)=F'_{j}(-\infty)=0.
\]
If, in addition, $\int_{-\infty}^\infty |t^{m-1}\phi_{j}(t)|dt<\infty$ for $2\le m\le 4$, then the Dominated Convergence theorem implies that $F'(x)$  is differentiable  and  thus for $1\le m\le 4$
\be \label{inv3}
F^{(m)}_{j}(x)=\frac{(-i)^{m-1}}{2\pi}\int_{-\infty}^\infty e^{-itx}t^{m-1}\phi_{k}(t)dt, \quad F^{(m)}_{j}(\infty)=F^{(m)}_{j}(-\infty)=0.
\ee

 Finally, we present some distributional relations for logistic order statistics with exponential shifts, which might be of independent interest. For similar identities see George and Rousseau \cite{GR87} (cf. Kotz et al. \cite{KKP01}, p.127).

\vspace{0.5cm}{\bf Lemma 1}

(i) Let $k$ and $m$ be two integers, such that $1\le k<m\le n$. Then
\be \label{lemma}
X_{m,n}-\sum_{j=k}^{m-1}\frac{E'_j}{j}\stackrel{d}{=}X_{k,n}+ \sum_{j=k}^{m-1}\frac{E''_j}{n-j}.
\ee
(ii) Let $1\le k\le n$. Then
\be \label{lemmaB}
X\stackrel{d}{=} X_{k,n} +\sum_{j=1}^{n-k}\frac{E'_j}{j}-\sum_{j=1}^{k-1}\frac{E''_j}{j},
\ee
where as usual $\sum_{j=1}^0 (\cdot)=0$.

{\bf Proof}.
Referring to (\ref{os_density9}),
 for the logistic c.d.f. $F_L(x)$ in (\ref{log_cdf}) we have
\beq \label{rec}
\phi_{k}(t)
        & = &
    k{n \choose k}\int_{-\infty}^\infty e^{itx}\left(\frac{e^x}{1+e^x}\right)^{k-1}\left(\frac{1}{1+e^{x}}\right)^{n-k}\frac{e^{-x}}{\left(1+e^{-x}\right)^2}dx   \nonumber \\
    & = &
    k{n \choose k}\int_0^1 u^{k-1+it}(1-u)^{n-k-it}du  \nonumber \\
    & = &
\frac{\Gamma(k+it)}{\Gamma(k)}\ \frac{\Gamma(n-k+1-it)}{\Gamma(n-k+1)}. \nonumber
\eeq
Denote the logistic ch.f. by
$\phi(t)=\Gamma(1+it)\Gamma(1-it)$.
Thus, the ch.f. of the logistic $k^{th}$  order statistic
for $2\le k\le n-1$ can be factored as follows
\beq \label{cf_r}
 \phi_{k}(t)   & = &
 \frac{(k-1+it)\ldots (1+it)\Gamma(1+it)}{(k-1)!}\ \frac{(n-k-it)\ldots (1-it)\Gamma(1-it)}{(n-k)!}  \nonumber \\
    & = &
    \left(1+\frac{it}{k-1}\right)\ldots(1+it)\left(1-\frac{it}{n-k}\right)\ldots (1-it)\phi(t).
\eeq
Similarly, referring to (\ref{os_density9}) with $k=1$ and $k=n$, we obtain
\be \label{cf_1n}
 \phi_{1}(t) = \prod_{j=1}^{n-1}\left(1-\frac{it}{j}\right)\phi(t)\quad \mbox{and}\quad
 \phi_{n}(t) = \prod_{j=1}^{n-1}\left(1+\frac{it}{j}\right)\phi(t).
\ee
(i) The relations (\ref{cf_r}) and (\ref{cf_1n}) imply for $1\le k<m\le n$
\nbeq
 \frac{\phi_{m}(t)}{\phi_{k}(t)}
 & = &
    \frac{\left(1-\frac{it}{n-k}\right)^{-1}\ldots \left(1-\frac{it}{n-m+1}\right)^{-1}}{\left(1+\frac{it}{m-1}\right)^{-1}\ldots \left(1+\frac{it}{k}\right)^{-1}}  \\
    & = &
    \frac{\phi_E\left(\frac{t}{n-k}\right)\ldots \phi_E\left(\frac{t}{n-m+1}\right)}{\phi_E\left(-\frac{t}{m-1}\right)\ldots \phi_E\left(-\frac{t}{k}\right)}, \nonumber
    \neeq
where $\phi_E(t)=(1-it)^{-1}$ is the standard exponential ch.f. This is equivalent to
\be   \label{cfs}
\phi_{m}(t)\prod_{j=k}^{m-1}\phi_E\left(-\frac{t}{j}\right)=
\phi_{k}(t)\prod_{j=k}^{m-1}\phi_E\left(\frac{t}{n-j}\right).
\ee
Since for independent random variables the product of their ch.f.'s corresponds to their sum, (\ref{cfs}) yields (\ref{lemma}). The assertion in (ii) follows directly from (\ref{cf_r}) and (\ref{cf_1n}).

\section{Adjacent Order Statistics: $r=1$}

To prove the sufficiency in this case
assume (\ref{lemmaA}) with $r=1$ and set $k_1=k$. (\ref{lemmaA}) yields
\be \label{cfs1}
\phi_{k+1}(t)\left(1-\frac{it}{n-k}\right) = \phi_{k}(t)\left(1+\frac{it}{k}\right).
\ee
Equation (\ref{cfs1}) and the inversion formula (\ref{inv3}) imply
\beq \label{cfs29}
F'_{k+1}(x)+\frac{1}{n-k}F''_{k+1}(x)
    & = &
  \frac{1}{2\pi}\int_{-\infty}^\infty e^{-itx}\left(1-\frac{it}{n-k}\right)\phi_{k+1}(t)dt \nonumber \\
  & = &
 \frac{1}{2\pi}\int_{-\infty}^\infty e^{-itx}\left(1+\frac{it}{k}\right)\phi_{k}(t)dt \nonumber \\
  & = &
   F'_{k}(x)-\frac{1}{k}F''_{k}(x).
   \eeq
Integrating (\ref{cfs29}), taking into account the boundary conditions, we obtain
\be \label{main1}
F_{k}(x)-F_{k+1}(x)=\frac{F'_{k}(x)}{k}+\frac{F'_{k+1}(x)}{n-k}.
\ee
Using (\ref{os_density9}) and (\ref{os_density3}), equation (\ref{main1}) can be simplified to
\[
F'(x)-F(x)(1-F(x))=0,
\]
i.e., (\ref{property}) holds, which along with the condition $F(0)=1/2$, completes the proof of the sufficiency.
The necessity follows directly from Lemma 1(i).


\section{Non-Adjacent Cases: $r=2$ and $r=3$}

Let us introduce the notation
\be \label{w19}
w(x)=\frac{F'(x)}{F(x)(1-F(x))},\qquad  x\in (l_F,r_F).
\ee
For notational simplicity, we will suppress the dependence of $w$, $F$, and $F_k$ on $x$.

\vspace{0.5cm}\paragraph{Proof of Theorem 1 when $r=2$.}
To start the proof of sufficiency, assume (\ref{lemmaA}) with $r=2$ and set $k_1=k$. Hence,
\be \label{cfs39}
\phi_{k+2}(t)\left(1-\frac{it}{n-k-1}\right)\left(1-\frac{it}{n-k}\right) = \phi_{k}(t)\left(1+\frac{it}{k}\right)\left(1+\frac{it}{k+1}\right).
\ee
Similarly to (\ref{cfs29}), using the inversion formula (\ref{inv3}),  one can see that (\ref{cfs39}) implies
\[
F'_{k}-F'_{k+2}=\sum_{j=k}^{k+1}\left(\frac{F''_k}{j}+\frac{F''_{k+2}}{n-j}\right)-\frac{F'''_k}{k(k+1)}+\frac{F'''_{k+2}}{(n-k)(n-k-1)}.
\]
Integrating, taking into account the boundary conditions, we obtain
\be \label{main1_m9}
F_{k}-F_{k+2}=\sum_{j=k}^{k+1}\left(\frac{F'_k}{j}+\frac{F'_{k+2}}{n-j}\right)-\frac{F''_k}{k(k+1)}+\frac{F''_{k+2}}{(n-k)(n-k-1)}.
\ee
Applying (\ref{os_density3}) to the left-hand side of (\ref{main1_m9}) yields
 \be \label{main11_m9}
 \sum_{j=k}^{k+1}{n \choose j}F^{j}(1-F)^{n-j}
   =  \sum_{j=k}^{k+1}\left(\frac{F'_k}{j}+\frac{F'_{k+2}}{n-j}\right)-\frac{F''_k}{k(k+1)}+\frac{F''_{k+2}}{(n-k)(n-k-1)}.
\ee
Using (\ref{os_density9}), we write the right-hand side of (\ref{main11_m9}) in terms of $F$, $F'$, and $F''$ only, and  after some algebra, obtain
\be \label{w39}
(2F-1)w'=P(n,F,w)(w-1)-k(2F-1)(w-1)^2,
\ee
where $w=w(x)$ is defined by (\ref{w19}) and $P(n,F,w)=(2nF-n-1)Fw-(n-1)F-1$.
Clearly $w\equiv 1$ is a solution of (\ref{w39}). This solution corresponds to the standard logistic c.d.f. $F_L$. To prove the sufficiency part of the theorem, it remains to show that this solution is unique.
According to the assumptions of the theorem, there are two distinct integers $k_1$ and $k_2$ such that (\ref{w39}) holds for both $k=k_1$ and $k=k_2$.
Writing (\ref{w39}) for these two values of $k$, subtracting the two equations from each other, and dividing by $k_1-k_2\ne 0$, we obtain
\be \label{final}
(2F-1)(w-1)^2=0.
\ee
Since $F(x)\neq 1/2$ for $x\ne 0$, the only solution of (\ref{final}) is $w(x)\equiv 1$, which implies the characteristic property (\ref{property}).
This, along with $F(0)=1/2$, completes the proof of the sufficiency. The necessity follows directly from Lemma~1(i).

\vspace{0.5cm}\paragraph{Proof of Theorem 1 when $r=3$.}
First, we will proof the sufficiency.
Assuming (\ref{lemmaA}) with $r=3$ and setting $k_1=k$, we have
\be \label{cfs3}
\phi_{k+3}(t)\prod_{j=0}^2\left(1-\frac{it}{n-k-j}\right)=
 \phi_{k}(t) \prod_{j=0}^2\left(1+\frac{it}{k+j}\right)
\ee
Similarly to (\ref{main1_m9}), using the inversion formula (\ref{inv3}) and integrating the resulting expression, one can see that (\ref{cfs3}) implies
\beq \label{main1_m}
F_{k}- F_{k+3} & = &
 \sum_{j=k}^{k+2}\left(\frac{F'_k}{j}+\frac{F'_{k+3}}{n-j}\right)-\sum\left(\frac{F''_k}{ij}-\frac{F''_{k+3}}{(n-i)(n-j)}\right)\\
 & & +\frac{F'''_k}{k(k+1)(k+2)}+\frac{F'''_{k+3}}{(n-k)(n-k-1)(n-k-2)},  \nonumber
\eeq
where the summation in $\sum$ is over all $k\le i<j\le k+2$.
Now, applying (\ref{os_density3}) to the left-hand side of (\ref{main1_m}) yields
\beq \label{main11_m}
  \sum_{j=k}^{k+2}{n \choose j}F^{j}(1-F)^{n-j}
 & = &
  \sum_{j=k}^{k+2}\left(\frac{F'_k}{j}+\frac{F'_{k+3}}{n-j}\right)-\sum\left(\frac{F''_k}{ij}-\frac{F''_{k+3}}{(n-i)(n-j)}\right) \\
 & & +\frac{F'''_k}{k(k+1)(k+2)}+\frac{F'''_{k+3}}{(n-k)(n-k-1)(n-k-2)},  \nonumber
\eeq
where the summation in $\sum$ is over all $k\le i<j\le k+2$.
Using (\ref{os_density9}), we write
$F^{(i)}_{k}$ and $F^{(i)}_{k+3}$ for $1\le i\le 3$ in terms of $F^{(i)}$ for $0\le i\le 3$. Then, after some tedious algebra, one can see that (\ref{main11_m}) is equivalent to
\be \label{w3}
Q(n,F,w)
  =  kP_1(F,w)(w-1)+kP_2(n,F,w)(w-1)^2+k^2P_3(F)(w-1)^3,
\ee
where $w=w(x)$ is defined by (\ref{w19}) and the polynomials  $P_1(F,w)=3(3F^2-3F+1)w'$, $P_2(n,F,w)=3(n-2)F^2-(n-8)F-3-2F[3nF^2-3(n+1)F+n+1]w$, $P_3(F)=(3F^2-3F+1)$, and $Q(n,F,w)$ do not depend on $k$.
One solution of (\ref{w3}) is $w\equiv 1$, which yields $F_L(x)$. It remains to prove that this is the only solution.
By the assumptions, (\ref{w3}) holds for three distinct values of $k$.
Writing (\ref{w3}) for $k=k_1$ and $k=k_2$ and  subtracting the two equations from each other, 
we obtain
\be \label{diff1}
P_1(F,w)(w-1)+P_2(n,F,w)(w-1)^2+(k_1+k_2)P_3(F)(w-1)^3=0.
\ee
Next, writing (\ref{w3}) for $k=k_2$ and $k=k_3$ and  subtracting the two equations from each other, 
we have
\be \label{diff2}
P_1(F,w)(w-1)+P_2(n,F,w)(w-1)^2+(k_2+k_3)P_3(F)(w-1)^3=0.
\ee
Finally, subtracting (\ref{diff1}) from (\ref{diff2}) and dividing by $k_3-k_1>0$,
we have
\be \label{diff3}
(3F^2-3F+1)(w-1)^3=0.
\ee
Since, $3F^2-3F+1\ne 0$, the only solution of (\ref{diff3}) is $w\equiv 1$, i.e., (\ref{property}) holds.
 Now, taking into account $F(0)=1/2$ we complete the proof of the sufficiency.
The necessity is a straightforward corollary of Lemma~1(i).

\section{Proof of Theorem 2}
\label{sec:5}

 In proving Theorem 2, we follow  the general scheme of the proof of Theorem 6 in Lin and Hu \cite{LH08}. Let, for the reader's convenience, recall below the notion of intensively monotone operators and strongly $\mathcal{E}$-positive families (cf. Kakosyan et al. \cite{KKM84}).

\vspace{0.5cm}{\bf Definition 1}. Let $\mathcal{C}=\mathcal{C}$ $[0,\infty)$ be the space of all real-valued functions defined and continuous in the interval $[0,\infty)$. The notation $f\ge g$ for $f,g\in \mathcal{C}$ means that $f(t)\ge g(t)$ for all $t\in [0,\infty)$. Let $A$ be an operator mapping some set $\mathcal{E}\subset \mathcal{C}$ into $\mathcal{C}$. We say that the operator $A$ is intensively monotone, if for any $f_1$ and $f_2$ belonging to $\mathcal E$, the condition $f_1(\tau)\ge f_2(\tau)$ for all $\tau\in (0,t)$ implies that $(Af_1)(\tau)\ge (Af_2)(\tau)$ for $\tau\in (0,t)$ and, in addition, the condition $f_1(\tau)>f_2(\tau)$ for all $\tau\in (0,t)$ implies that $(Af_1)(t)>(Af_2)(t)$.

\vspace{0.5cm}{\bf Definition 2}. Let $\mathcal{E}\subset \mathcal{C}$ and $\{f_\lambda\}_{\lambda\in \Lambda}$ be a family of elements of $\mathcal{E}$. We say that the family $\{f_\lambda\}_{\lambda\in \Lambda}$ is strongly $\mathcal{E}$-positive if the following conditions hold:

(i) for any $f\in \mathcal{E}$ there are $t_0>0$ and $\lambda_0\in \Lambda$ such that $f(t_0)=f_{\lambda_0}(t_0)$;

(ii) for any $f\in \mathcal{E}$ and any $\lambda\in \Lambda$ either $f(t)=f_\lambda(t)$ for all $t\in [0,\infty)$, or there is  $\delta>0$ such that the difference $f(t)-f_\lambda(t)$ does not vanish (preserves its sign) in the interval $(0,\delta]$.

\vspace{0.5cm}The following example and lemma play essential role in the proof of Theorem~2.

\vspace{0.5cm}{\bf Example 1}. (Lin and Hu \cite{LH08}) Define $\mathcal{E}$ to be the set of all survivor functions $\overline{G}$, which are real, analytic, and strictly decreasing in $[0,\infty)$. Assume also that for each $k$, the $k$th derivative $\overline{G}^{(k)}$ is strictly monotone in some interval $[0, \delta_k)$. Denote $\overline{G}_\lambda(x)=1/(1+\lambda x)$, $x>0$, where $\lambda\in \Lambda=(0,\infty)$. It is proved in Lemma 5 of Lin and Hu (2008) that the family $\{\overline{G}_\lambda\}_{\lambda\in \Lambda}$ is strongly $\mathcal{E}$-positive.

\vspace{0.5cm}{\bf Lemma 2} (Kakosyan et al. \cite{KKM84}, Theorem 1.1.1). Let $A$ be an intensively monotone operator on $\mathcal{E}\subset\mathcal{C}$ and let $\{f_\lambda\}_{\lambda\in \Lambda}$ be a strongly $\mathcal{E}$-positive family. Assume that $Af_\lambda=f_\lambda$ for all $\lambda\in \Lambda$. Then the condition $Af=f$, where $f\in \mathcal{E}$, implies that there is $\lambda\in \Lambda$ such that $f=f_\lambda$. In other words, all solutions of the equation $Af=f$, belonging to $\mathcal{E}$, coincide with elements of the family $\{f_\lambda\}_{\lambda\in \Lambda}$.

\vspace{0.5cm}\paragraph{Proof of Theorem 2} The case $k=1$ is given as Theorem 6 in Lin and Hu \cite{LH08}. We assume that $2\le k\le n$. The necessity statement follows from Lemma~1(ii). Next we will prove the sufficiency.
Let, as before, $E'_j$ and  $E''_j$ for $1\le j\le n$ be i.i.d. standard exponential.  Recall  the well-known (e.g., David and Nagaraja (2003), pp.17-18) formula for the maximum $E'_{n,n}:=\max\{E'_1, \ldots E'_n\}$,
\be \label{nn_ident}
E'_{n,n}\stackrel{d}{=}\sum_{j=1}^n \frac{E'_j}{j}, \qquad n\ge 1.
\ee
Note also that
$U'_j :=\exp\{-E'_j\}$ for $1\le j\le n$ are uniform on $[0,1]$. Therefore, $U'_{1,n}:=\min\{U'_1,\ldots, U'_n\}\stackrel{d}{=}\exp\{-E'_{n,n}\}$ is distributed as the minimum of $n$ i.i.d.  uniform $[0,1]$ variables. Now, (\ref{goal1}) and (\ref{nn_ident}) yield
\beq \label{V_ident}
\exp\{X\} & \stackrel{d}{=} & \frac{\exp\{X_{k,n}\}\exp\{ -\sum_{j=1}^{k-1}E''_j/j\}}{\exp\{ -\sum_{j=1}^{n-k}E'_j/j\}}\\
    & \stackrel{d}{=} &
    \frac{\exp\{X_{k,n}\}\exp\{ -E''_{k-1,k-1}\}}{\exp\{ -E'_{n-k,n-k}\}} \nonumber \\
    & \stackrel{d}{=} &
    \frac{U''_{1,k-1}\exp\{X_{k,n}\}}{U'_{1,n-k}}. \nonumber
\eeq
 Let us make the change of variables
\[
\xi=\frac{U''_{1,k-1}\exp\{X_{k,n}\}}{U'_{1,n-k}},\qquad \eta=U''_{1,k-1}, \qquad \zeta=U'_{1,n-k}.
\]
 Denote $G_X(x):=P\left(\exp\{X\}\le x\right)$ and $\overline{G}_X(x):=1-G_X(x)$. For the p.d.f. of $\exp\{X_{k,n}\}$ and $U_{1,j}$, we have
\[
f_{\exp\{X_{k,n}\}}(x)=k{n \choose k}G^{k-1}_X(x)\overline{G}^{n-k}_X(x)G'_X(x),\quad f_{U_{1,j}}(y)=j(1-y)^{j-1}.
\]
Therefore, for the density $f_\xi(u)$ of $\xi$, we obtain
\beq \label{xi}
f_\xi(u)& = &
    \int_0^1\int_0^1 f_{\xi,\eta,\zeta}(u,v,w)dv dw \\
      & = &
     \int_0^1\int_0^1 C(v,w)\ G^{k-1}_X\left(\frac{uw}{v}\right)\overline{G}^{n-k}_X\left(\frac{uw}{v}\right)G'_X\left(\frac{uw}{v}\right)\frac{w}{v}dv dw, \nonumber
  \eeq
  where $C(v,w)=n!(1-v)^{k-2}(1-w)^{n-k-1}/(k-2)!(n-k-1)!$.
Now, referring to (\ref{V_ident}) and (\ref{xi}), one can see that $\overline{G}_X(x)$ satisfies for $n\ge 3$
\beq \label{k=2}
\overline{G}_X(x) & = & \int_x^\infty f_\xi(u)du\\
& = &
\int_0^1 \int_0^1 \frac{C(v,w)}{k}\ \left[ \int_{G_X(xw/v)}^\infty  \overline{G}^{n-k}_X\left(\frac{uw}{v}\right)dG^k_X\left(\frac{uw}{v}\right)\right] dv dw, \nonumber \\
& = &
\int_0^1 \int_0^1 \frac{C(v,w)}{k}\  H_k\left(\overline{G}_X\left(\frac{xw}{v}\right)\right)dvdw, \nonumber
    \eeq
where, making the substitution $t=G_X(uw/v)$, one can write $H_k(y)$ as
\be \label{H_k}
H_k(y)=\int_{1-y}^1 (1-t)^{n-k}dt^k,\qquad 0<y<1.
\ee
Define $\mathcal{E}$ to be the set of all survival functions $\overline{G}$ as in Example 1 above. Define also an operator $A$ on $\mathcal{E}$ by
\be \label{oper}
A\overline{G}_X(x)
=
\int_0^1 \int_0^1 C(v,w)\ H_k\left(\overline{G}_X\left(\frac{xw}{v}\right)\right)dvdw
\ee
Note that
$H'_k(y)=ky^{n-k}(1-y)^{k-1}>0$ for $0<y<1$ and hence $H_k(y)$ is an increasing function  for any $y\in (0,1)$ and $1\le k\le n$. Now, one can see that $A$
is an intensively monotone operator on $\mathcal{E}$. Furthermore, by Example~1, the family $\{\overline{G}_\lambda\}_{\lambda\in \Lambda}$ is strongly $\mathcal{E}$-positive. We have that $A\overline{G}_\lambda=\overline{G}_\lambda$ by the necessity part of the theorem. Finally, by Lemma 2 and (\ref{k=2}) (which means $A\overline{G}_X=\overline{G}_X$), we conclude that $\overline{G}_X=\overline{G}_\lambda$. The condition $F(0)=1/2$ implies $\lambda=1$ and completes the proof.

\section{Concluding Remarks}

Theorem 1 is a characterization result for standard logistic distribution based on distributional equalities between two order statistics plus or minus sums of independent exponential variables (two-sided exponential shifts). In case of adjacent order statistics, the characterization involves a single equality. If the order statistics are two or three spacings away - two or three equalities are needed, respectively. Our second result, characterizes the logistic distribution by certain decomposition of the parent variable $X$ into one order statistic plus a linear combination of independent exponential variables. In the corollary we singled out the particular case of the sample median.

One open question is if the number of characterizing equalities in the non-adjacent cases of Theorem~1 can be reduced. Another area of future work is to study other identities involving order statistics subject to exponential shifts. Some results in this direction are given in Ahsanullah et al. \cite{ANY11}. Finally, one can look for connections between random shifts and contractions, see for example Weso{\l}owski and  Ahsanullah (2004).

\end{document}